\documentclass[11pt,reqno]{amsart}
\usepackage{amsfonts}
\usepackage{mathrsfs}
\usepackage{amsmath}
\usepackage{amssymb,amsfonts}

\newtheorem{thm}{Theorem}[section]
\newtheorem{lem}{Lemma}[section]

\numberwithin{equation}{section}

\def\pf{{\textit {Proof:} }}

\newcommand{\mysection}[1]{\section{#1}\setcounter{equation}{0}}

\newfont{\bb}{msbm10 at 11pt}

\newcommand{\bal}{\begin{aligned}}      \newcommand{\eal}{\end{aligned}}
\newcommand{\ba}{\begin{array}}      \newcommand{\ea}{\end{array}}
\newcommand{\bc}{\begin{center}}     \newcommand{\ec}{\end{center}}
\newcommand{\be}{\begin{enumerate}}  \newcommand{\ee}{\end{enumerate}}
\newcommand{\beq}{\begin{eqnarray}}  \newcommand{\eeq}{\end{eqnarray}}
\newcommand{\beQ}{\begin{eqnarray*}} \newcommand{\eeQ}{\end{eqnarray*}}
\newcommand{\bi}{\begin{itemize}}    \newcommand{\ei}{\end{itemize}}
\newcommand{\bt}{\begin{tabular}}    \newcommand{\et}{\end{tabular}}
\newcommand{\bdm}{\begin{displaymath}} \newcommand{\edm}{\end{displaymath}}

\newcommand{\rw}{\rightarrow}

\newcommand{\ep}{\epsilon}

\def\qed{\hfill{Q.E.D.}\smallskip}

\newcommand{\ls}{\setlength{\baselineskip}{12pt}
                 \setlength{\parskip}{3mm}}

\begin{document}

\title[]{Geodesics on metrics of self-dual Taub-Nut type}
\author[C Liu]{Chuxiao Liu$^{\flat}$}
\address[]{$^{\flat}$School of Mathematics and Information Science, Guangxi University, Nanning, Guangxi 530004, PR China}
\address[]{$^{\flat}$Guangxi Center for Mathematical Research, Guangxi University, Nanning, Guangxi 530004, PR China}
\address[]{$^{\flat}$Guangxi Base, Tianyuan Mathematical Center in Southwest China, Nanning, Guangxi 530004, PR China}
\email{cxliu@gxu.edu.cn}
\author[Q Pu]{Qingtao Pu$^{\dag}$}
\address[]{$^{\dag}$School of Mathematics and Information Science, Guangxi University, Nanning, Guangxi 530004, PR China}
\email{2206301035@gxu.edu.cn}
\date{}

\begin{abstract}
Geodesic equations are solved when at least two
of $\tau$, $\theta$, $\varphi$ are constant on
metrics of self-dual Taub-NUT type. They can also be solved
also on self-dual Taub-NUT metrics if only
$r$, $\theta$ or $\varphi$ is constant. However, the explicit solution of the geodesic
equations is not available yet if only $\tau$ is constant.
\end{abstract}

\maketitle \pagenumbering{arabic}

\mysection{Introduction}\ls
Geodesics are fundamental geometric objects in differential geometry and play a crucial role in the study of general relativity. They describe the motion of light and massive particles through null-like and time-like geodesics, respectively. Geodesics have been extensively studied in various specific spacetimes (see, e.g., \cite{C,CM} and references therein). Recently, Battista and Esposito analyzed geodesic motion in Euclidean Schwarzschild geometry \cite{BE}, while Yang and Zhang explored the geodesics of metrics of Eguchi-Hanson type \cite{YZ}. Both metrics can be viewed as generalizations of the Schwarzschild metric and belong to the class of gravitational instantons. Gravitational instantons are defined as nonsingular, complete, positive-definite solutions of the classical vacuum Einstein equations or the Einstein equations with a cosmological constant term, playing a significant role in the Euclidean approach to quantum gravity \cite{GH}. Therefore, studying geodesic motion on other gravitational instantons is of great interest. This paper aims to investigate the geodesic motion on the self-dual Taub-NUT metric, which also belongs to the class of gravitational instantons \cite{GH}.

The Taub-NUT metric is a solution to Einstein's field equations in general relativity that extends the Taub solution \cite{T} and incorporates elements of the Newman-Unti-Tamburino (NUT) metric \cite{NUT}. It describes a four-dimensional spacetime characterized by mass and a parameter known as the NUT charge, linked to a ``gravitational magnetic monopole'' or twisting geometry \cite{GH,NUT,T}. This metric is notable for its exotic features, including closed timelike curves (suggesting the theoretical possibility of time travel) and the Misner string, a topological defect similar to the Dirac string in electromagnetism \cite{CG,CGG}. The Taub-NUT spacetime plays a pivotal role in the study of gravitational instantons, black hole thermodynamics, and string theory, with numerous related results available (see, e.g., \cite{LL} and references therein). The geodesic completeness of this metric has already been demonstrated in \cite{CG,CGG}, motivating us to solve the geodesic equations through the removable singularity.

The paper is organized as follows. In Sect. 2, we give a quick review of the self-dual Taub-NUT metrics. We also show the geodesic equaitons of metrics of self-dual Taub-NUT type. In Sect. 3, we solve the geodesic equations of metrics of self-dual Taub-NUT type in following cases: (1) Sect. 3.1, where $\tau$,$\theta$,$\varphi$ are constants, geodesics pass through the removable singularity; (2) Sect. 3.2, where $\theta$,$\varphi$ are constants, geodesics do not pass through the removable singularity. If geodesics pass through the removable singularity, $\tau$ must be constant and the geodesic equations reduce to Sect. 3.1; (3) Sect. 3.3, where $\theta\in (0,\pi), \tau $are constants. It yields that $\theta=\frac{\pi}{2}$ if $\varphi$ and $r$ are not constants; (4) Sect. 3.4, where $\tau, \varphi$ are constants, geodesics do not pass through the removable singularity. If geodesics pass through the removable singularity, $\theta$ must be constant; (5) Sect. 3.5, where $r$ is constant, geodesics pass through the removable singularity; (6) Sect. 3.6, where $\varphi$ is constant, the geodesic equations reduce to Sect. 3.2 or Sect. 3.3; (7) Sect. 3.7, where $\theta\in(0,\pi)$ is constant, geodesics do not pass through the removable singularity. If geodesics pass through the removable singularity, $\varphi$ must be constant and the geodesic equations reduce to Sect. 3.2. We point out that the solution of geodesics equations is not available if only $\tau$ is constant.

\mysection{Metrics of self-dual Taub-NUT type and geodesic equations}\ls

In this section,we introduce the metrics of self-dual Taub-NUT type and show their geodesic equations. Metrics of self-dual Taub-NUT type are given by \cite{GH}
\beQ
\begin{aligned}
ds^2=\frac{r-n}{r+n}(d\tau+2n\cos\theta d\varphi)^2+(r^2-n^2)(d\theta^2+\sin^2\theta d\varphi^2)+\frac{r+n}{r-n}dr^2.
\end{aligned}
\eeQ
The apparent singularity at $r=n$ is removable and just the origin in hyperspherical polar coordinates. The topology of the manifold is $R^4$. And the surfaces of constant $r>n$ are equipped with the topology of 3-spheres with the Euler angles $(\tau(2n)^{-1},\theta,\phi)$. The curvature is self-dual with the following positively oriented orthonormal basis \cite{GH}
\begin{align*}
\omega^0&=\sqrt{\frac{r+n}{r-n}}dr,\\
\omega_1&=\sqrt{r^2-n^2}\left(\cos\frac{\tau}{2n}d\theta+\sin\frac{\tau}{2n}\sin\theta d\phi\right),\\
\omega^2&=\sqrt{r^2-n^2}\left(-\sin\frac{\tau}{2n}d\theta+\cos\frac{\tau}{2n}\sin\theta d\phi\right),\\
\omega^3&=\sqrt{\frac{r-n}{r+n}}(d\tau+2n\cos\theta d\phi).
\end{align*}
The nontrivial metric components are
\begin{align*}
g_{\tau\tau}&=\frac{r-n}{r+n},\,\,g_{\tau\varphi}=g_{\varphi\tau}=\frac{2n\cos\theta(r-n)}{r+n},\\
g_{\varphi\varphi}&=\frac{4n^2\cos^2\theta(r-n)}{r+n}+(r^2-n^2)\sin^2\theta,\\
g_{rr}&=\frac{r+n}{r-n},\,\, g_{\theta\theta}=r^2-n^2,\\
g^{\tau\tau}&=\frac{4n^2\cos^2\theta}{(r^2-n^2)\sin^2\theta}+\frac{r+n}{r-n},\\
g^{\tau\varphi}&=g^{\varphi\tau}=\frac{-2n\cos\theta}{(r^2-n^2)\sin^2\theta},\,\, g^{\varphi\varphi}=\frac{1}{(r^2-n^2)\sin^2\theta},\\
g^{rr}&=\frac{r-n}{r+n},\,\, g^{\theta\theta}=\frac{1}{r^2-n^2}.
\end{align*}
The nontrivial Christoffel symbols are
\beQ
\begin{aligned}
\Gamma_{\tau r}^\tau &=\frac{n}{r^2-n^2},\Gamma_{\tau\theta}^\tau =\frac{2n^2\cos\theta}{(r+n)^2\sin\theta},\Gamma_{\varphi r}^\tau=-\frac{2n\cos\theta}{r+n},\\
\Gamma_{\varphi\theta}^\tau &=\frac{4n^3\cos^2\theta-n\sin^2\theta(r+n)^2-2n(r+n)^2\cos^2\theta}{(r+n)^2\sin\theta},\\
\Gamma_{\tau\tau}^r &=-\frac{n(r-n)}{(r+n)^3},\Gamma_{\tau\varphi}^r =-\frac{2n^2(r-n)\cos\theta}{(r+n)^3},\Gamma_{rr}^r =-\frac{n}{r^2-n^2},\\
\Gamma_{\theta\theta}^r &=-\frac{r(r-n)}{r+n},\Gamma_{\varphi\varphi}^r =-\left[\frac{4n^3\cos^2\theta}{(r+n)^2}+r\sin^2\theta\right]\frac{r-n}{r+n},\\
\Gamma_{\tau\varphi}^\theta &=\frac{n\sin\theta}{(r+n)^2},\Gamma_{r \theta}^\theta =\frac{r}{r^2-n^2},\\
\Gamma_{\varphi\varphi}^\theta &=\frac{4n^2\cos\theta\sin\theta}{(r+n)^2}-\sin\theta\cos\theta,\\
\Gamma_{\tau\theta}^\varphi &=-\frac{n}{(r+n)^2\sin\theta},\Gamma_{\varphi r}^\varphi =\frac{r}{r^2-n^2},\\
\Gamma_{\varphi\theta}^\varphi &=-\frac{2n^2\cos\theta}{(r+n)^2\sin\theta}+\frac{\cos\theta}{\sin\theta}.\\
\end{aligned}
\eeQ
Thus geodesic equations for parameter t read as
\beQ
\begin{aligned}
\frac{d^2\tau}{dt^2}&+\frac{4n^2\cos\theta}{(r+n)^2\sin\theta}\frac{d\tau}{dt}\frac{d\theta}{dt}+\frac{2n}{r^2-n^2}\frac{d\tau}{dt}\frac{dr}{dt}\\&
+\frac{2[4n^3\cos^2\theta-n\sin^2\theta(r+n)^2-2n(r+n)^2\cos^2\theta]}
{(r+n)^2\sin\theta}\frac{d\varphi}{dt}\frac{d\theta}{dt}\\
&-\frac{4n\cos\theta}{r+n}\frac{d\varphi}{dt}\frac{dr}{dt}=0,
\end{aligned}
\eeQ
\beq
\begin{aligned}
\frac{d^2\varphi}{dt^2}&+2\left[-\frac{2n^2\cos\theta}{(r+n)^2\sin\theta}+\frac{\cos\theta}{\sin\theta}\right]\frac{d\varphi}{dt}\frac{d\theta}{dt}\\
&+\frac{2r}{r^2-n^2}\frac{d\varphi}{dt}\frac{dr}{dt}-\frac{2n}{(r+n)^2\sin\theta}\frac{d\tau}{dt}\frac{d\theta}{dt}=0,\\
\frac{d^2\theta}{dt^2}&+\frac{2r}{r^2-n^2}\frac{dr}{dt}\frac{d\theta}{dt}+
\left[\frac{4n^2\cos\theta\sin\theta}{(r+n)^2}-\sin\theta\cos\theta\right]\left(\frac{d\varphi}{dt}\right)^2\\
&+\frac{2n\sin\theta}{(r+n)^2}\frac{d\tau}{dt}\frac{d\varphi}{dt}=0,\\
\frac{d^2r}{dt^2}&-\frac{n}{r^2-n^2}\left(\frac{dr}{dt}\right)^2-\frac{n(r-n)}{(r+n)^3}\left(\frac{d\tau}{dt}\right)^2-\frac{r(r-n)}{r+n}\left(\frac{d\theta}{dt}\right)^2\\
&-\left[\frac{4n^3\cos^2\theta}{(r+n)^2}+r\sin^2\theta\right]\frac{r-n}{r+n}\left(\frac{d\varphi}{dt}\right)^2-\frac{4n^2(r-n)\cos\theta}{(r+n)^3}\frac{d\tau}{dt}\frac{d\varphi}{dt}=0.\\
\end{aligned}\label{ge}
\eeq

Throughout the paper, we denote
\begin{align*}
\ep=\pm 1.
\end{align*}
\section{Geodesics on metrics of self-dual Taub-NUT type}

In this section, we shall solve the geodesic equations for the following cases.

\subsection{Geodesics for constant $\tau$,$\theta$,$\varphi$}\ls

In this case, the geodesic equations reduce to\\
\begin{equation}
\begin{aligned}
&\frac{d^2r}{dt^2}-\frac{n}{r^2-n^2}\left(\frac{dr}{dt}\right)^2=0.    \label{1}\\
\end{aligned}
\end{equation}
\begin{thm}\label{t1} The geodesics for metrics of self-dual Taub-NUT type with constant $\theta$,$\varphi$,$\tau$ and passing through $r=n$ with conditions
\beQ
\lim\limits_{r \rw n} {t}= t_{1},\quad \lim\limits_{r \rw n}{\sqrt{\frac{r+n}{r-n}}\frac{dr}{dt}}=r_1.
\eeQ
satisfy
\beQ
t(r)=t_1+\frac{1}{r_1}\left[\sqrt{r^2-n^2}+2n\ln {\left (\sqrt{r+n}+\sqrt{r-n}\right )}\right].
\eeQ
\end{thm}
\pf (\ref{1}) implies that
\beQ
\frac{d}{dt}\left(\sqrt{\frac{r+n}{r-n}}\frac{dr}{dt}\right)=0.
\eeQ
Therefore,
\beQ
\sqrt{\frac{r+n}{r-n}}\frac{dr}{dt}=r_1.
\eeQ
Thus,
\beQ
\frac{dt}{dr}=\frac{1}{r_1}\sqrt{\frac{r+n}{r-n}}.
\eeQ
The theorem follows by integrating it from $n$ to $r$.\qed

\subsection{Geodesics for constant $\theta$,$\varphi$}\ls
The geodesic equations reduce to
\begin{subequations}
\begin{align}
\frac{d^2r}{dt^2}&-\frac{n}{r^2-n^2}\left(\frac{dr}{dt}\right)^2-\frac{n(r-n)}{(r+n)^3}\left(\frac{d\tau}{dt}\right)^2=0, \label{2}\\
\frac{d^2\tau}{dt^2}&+\frac{2n}{r^2-n^2}\frac{d\tau}{dt}\frac{dr}{dt}=0.\label{3}
\end{align}
\end{subequations}

Let $r_1,\tau_0$ be constant and $r_1\neq 0$. Denote
\begin{align*}
R_1=n+\frac{2\tau_0^2n}{r_1^2}.
\end{align*}
Obviously, $R_1\geq n$, and $R_1=n$ if and only if $\tau_0=0$.
\begin{thm}\label{t2} Let $r\geq R_1>n$. The geodesics for metrics of self-dual Taub-NUT type with constant $\theta$,$\varphi$ and conditions
\beQ
\begin{aligned}
&\lim\limits_{r \rw R_1} t= t_{1},\,\,\lim\limits_{r \rw R_1} \tau = \tau_{1},\,\,\lim\limits_{r \rw R_1} {\frac{r-n}{r+n}\frac{d\tau}{dt}} = \tau_{0},\\
&\lim\limits_{r\rw R_1}{\left[\frac{r+n}{r-n}\left(\frac{dr}{dt}\right)^2+g_1\right]}=r_1^2.\\
\end{aligned}
\eeQ
satisfy
\beQ
\begin{aligned}
t(r)=&t_1+\frac{\ep}{r_1}\left[\sqrt{r-R_1}\sqrt{r+n}+(n+R_1)\ln \left(\sqrt{r-R_1}+\sqrt{r+n}\right)\right],\\
\tau(r)=&\tau_1+\frac{\ep \tau_0}{r_1}\Bigg[(5n+R_1)\ln (\sqrt{r-R_1}+\sqrt{r+n})\\
&+\frac{2(2n)^{\frac{3}{2}}}{\sqrt{R_1-n}}\arctan \left(\frac{\sqrt{2n}\sqrt{r-R_1}}{\sqrt{R_1-n}\sqrt{r+n}}\right)
+\sqrt{r-R_1}\sqrt{r+n}\Bigg].
\end{aligned}
\eeQ
\end{thm}
\pf The geodesic equation (\ref{3}) implies that
\begin{align*}
\frac{d}{dt}\left(\frac{r-n}{r+n}\frac{d\tau}{dt}\right)=0.
\end{align*}
Therefore,
\begin{equation*}
\begin{aligned}
&\frac{d\tau}{dt}=\tau_0\frac{r+n}{r-n}.
\end{aligned}
\end{equation*}
Substituting it into (\ref{2}), we obtain
\begin{equation}
\begin{aligned}
&\frac{d^2r}{dt^2}-\frac{n}{r^2-n^2}\left(\frac{dr}{dt}\right)^2-\frac{\tau_0^2 n}{r^2-n^2}=0.  \label{4}\\
\end{aligned}
\end{equation}
Similarly, let
\begin{align}
h_1=\frac{r+n}{r-n}.\label{h1}
\end{align}
Then
$$\frac{h_1'}{2h_1}=\frac{-n}{r^2-n^2}.$$
Thus, (\ref{4}) implies that
\begin{equation}
\begin{aligned}
&\frac{d}{dt}\left[h_1\left(\frac{dr}{dt}\right)^2+\frac{2\tau_0^2n}{r-n}\right]=0. \label{5}
\end{aligned}
\end{equation}
Thus,
\beQ
h_1\left(\frac{dr}{dt}\right)^2+\frac{2\tau_0^2n}{r-n}=r_1^2.
\eeQ
Therefore,
\beQ
\begin{aligned}
\frac{dt}{dr}&=\frac{\ep}{r_1}\sqrt{\frac{r+n}{r-R_1}},\\
\frac{d\tau}{dr}&=\frac{d\tau}{dt}\frac{dt}{dr}=\frac{\ep \tau_0}{r_1}\frac{(r+n)^{\frac{3}{2}}}{(r-n)\sqrt{r-R_1}}.
\end{aligned}
\eeQ
The theorem follows by integrating the above equations from $R_1$ to $r$. \qed

\subsection{Geodesics for constant $\theta\in (0,\pi)$, $\tau$}\ls
The geodesic equations reduce to
\begin{align}
&\frac{d^2r}{dt^2}-\frac{n}{r^2-n^2}\left(\frac{dr}{dt}\right)^2-\left[\frac{4n^3\cos^2\theta}{(r+n)^2}+r\sin^2\theta\right]\frac{r-n}{r+n}\left(\frac{d\varphi}{dt}\right)^2=0,\nonumber\\
&\frac{d^2\varphi}{dt^2}+\frac{2r}{r^2-n^2}\frac{d\varphi}{dt}\frac{dr}{dt}=0,\nonumber\\
&\frac{4n\cos\theta}{r+n}\frac{d\varphi}{dt}\frac{dr}{dt}=0,  \label{7}\\
&\left[\frac{4n^2\cos\theta\sin\theta}{(r+n)^2}-\sin\theta\cos\theta\right]\left(\frac{d\varphi}{dt}\right)^2=0. \label{8}
\end{align}
(\ref{7}) and (\ref{8}) imply that
\begin{align*}
\theta=\frac{\pi}{2},\,\, \text{or}\,\, \frac{d\varphi}{dt}=0,\,\, \text{or}\,\, \frac{dr}{dt}=0.
\end{align*}
It reduce to Theorem \ref{t1} if $\frac{d\varphi}{dt}=0$, while $r,\varphi$ are constants if $\frac{dr}{dt}=0$.

Now we focus on the case $\theta=\frac{\pi}{2}$. The geodesic equations reduce to
\begin{align}
\frac{d^2r}{dt^2}&-\frac{n}{r^2-n^2}\left(\frac{dr}{dt}\right)^2-\frac{r(r-n)}{r+n}\left(\frac{d\varphi}{dt}\right)^2=0,\label{9}\\
\frac{d^2\varphi}{dt^2}&+\frac{2r}{r^2-n^2}\frac{d\varphi}{dt}\frac{dr}{dt}=0. \label{10}
\end{align}
Let $r_1, \varphi_0$ be constant and $r_1\neq 0$. Denote
\begin{align*}
R_2=\sqrt{n^2+\frac{\varphi_0^2}{r_1^2}}.
\end{align*}
Obviously, $R_2\geq n$ and $R_2=n$ if and only if $\varphi_0=0$.
\begin{thm}\label{t3} The geodesics for metrics of self-dual Taub-NUT type with constant $\theta=\frac{\pi}{2}$,$\tau$ with conditions
\beQ
\begin{aligned}
&\lim\limits_{r \rw R_2} t= t_{1},\,\,\lim\limits_{r \rw R_2} \varphi = \varphi_{1},\,\,\lim\limits_{r \rw R_2} {(r^2-n^2)\frac{d\varphi}{dt}} = \varphi_{0},\\
&\lim\limits_{r\rw R_2}{\left[\frac{r+n}{r-n}\left(\frac{dr}{dt}\right)^2+\frac{\varphi_0^2}{r^2-n^2}\right]}=r_1^2.\end{aligned}
\eeQ
satisfy
\beQ
\begin{aligned}
t(r)&=t_1+\frac{\ep}{r_1}\left[\sqrt{r^2-R_2^2}+n\ln \left(r+\sqrt{r^2-R_2^2}\right)\right],\\
\varphi(r)&=\varphi_1+\ep\arctan \frac{r_1(rn-R_2^2)}{\varphi_0\sqrt{r^2-R_2^2}}.
\end{aligned}
\eeQ
\end{thm}
\pf (\ref{10}) implies that
\begin{align*}
\frac{d}{dt}\left[(r^2-n^2)\frac{d\varphi}{dt}\right]=0.
\end{align*}
Thus,
\begin{align*}
(r^2-n^2)\frac{d\varphi}{dt}=\varphi_0.
\end{align*}
Therefore,
\begin{align*}
\frac{d\varphi}{dt}=\frac{\varphi_0}{r^2-n^2}.
\end{align*}
Substituting it into (\ref{9}) ,we obtain
\begin{equation*}
\begin{aligned}
&\frac{d^2r}{dt^2}-\frac{n}{r^2-n^2}\left(\frac{dr}{dt}\right)^2-\frac{r\varphi_0^2}{(r-n)(r+n)^3}=0.
\end{aligned}
\end{equation*}
Then we have
\begin{align*}
\frac{d}{dt}\left[\frac{r+n}{r-n}\left(\frac{dr}{dt}\right)^2+\frac{\varphi_0^2}{r^2-n^2}\right]=0.
\end{align*}
Thus,
\begin{align*}
\left[\frac{r+n}{r-n}\left(\frac{dr}{dt}\right)^2+\frac{\varphi_0^2}{r^2-n^2}\right]=r_1^2.
\end{align*}
This implies that
\begin{align*}
\frac{dr}{dt}=\ep\frac{\sqrt{r_1^2(r^2-n^2)-\varphi_0^2}}{r+n}.
\end{align*}
Therefore,
\beQ
\begin{aligned}
\frac{dt}{dr}&=\frac{\ep}{r_1}\frac{r+n}{\sqrt{r^2-R_2^2}},\\
\frac{d\varphi}{dr}&=\frac{d\varphi}{dt}\frac{dt}{dr}=\frac{\ep\varphi_0}{r_1(r-n)\sqrt{r^2-R_2^2}}.
\end{aligned}
\eeQ
The theorem follows by integrating the above equations from $R_2$ to $r$.\qed

\subsection{Geodesics for constant $\varphi$,$\tau$}\ls
The geodesic equations reduce to
\begin{align}
\frac{d^2r}{dt^2}&-\frac{n}{r^2-n^2}\left(\frac{dr}{dt}\right)^2-\frac{r(r-n)}{r+n}\left(\frac{d\theta}{dt}\right)^2=0,\label{14}\\
\frac{d^2\theta}{dt^2}&+\frac{2r}{r^2-n^2}\frac{dr}{dt}\frac{d\theta}{dt}=0.\label{15}
\end{align}
Let $r_1, \theta_0$ be constant and $r_1\neq 0$. Denote
\begin{align*}
R_3=\sqrt{n^2+\frac{\theta_0^2}{r_1^2}}.
\end{align*}
Then $R_3\geq n$ and $R_3=n$ if and only if $\theta_0=0$.
\begin{thm}\label{t4} The geodesics for metrics of self-dual Taub-NUT type with constant $\varphi$,$\tau$ with conditions
\beQ
\begin{aligned}
&\lim\limits_{r \rw R_3} t= t_{1},\,\,\lim\limits_{r \rw R_3} \theta = \theta_{1},\,\,\lim\limits_{r \rw R_3} {(r^2-n^2)\frac{d\theta}{dt}} = \theta_{0},\\
&\lim\limits_{r\rw R_3}{\left[\frac{r+n}{r-n}\left(\frac{dr}{dt}\right)^2+\frac{\theta_0^2}{r^2-n^2}\right]}=r_1^2.
\end{aligned}
\eeQ
satisfy
\beQ
\begin{aligned}
t(r)&=t_1+\frac{\ep}{r_1}\left[\sqrt{r^2-R_3^2}+n\ln \left(r+\sqrt{r^2-R_3^2}\right)\right],\\
\theta(r)&=\theta_1+\ep\arctan \frac{r_1(rn-R_3^2)}{\theta_0\sqrt{r^2-R_3^2}}.
\end{aligned}
\eeQ
\end{thm}
\pf (\ref{15}) implies that
\begin{align*}
\frac{d}{dt}\left[(r^2-n^2)\frac{d\theta}{dt}\right]=0.
\end{align*}
This leads to
\begin{align*}
(r^2-n^2)\frac{d\theta}{dt}=\theta_0.
\end{align*}
Therefore,
\begin{align*}
\frac{d\theta}{dt}=\frac{\theta_0}{r^2-n^2}.
\end{align*}
Substituting it into (\ref{14}), we obtain
\begin{equation}
\begin{aligned}
&\frac{d^2r}{dt^2}-\frac{n}{r^2-n^2}\left(\frac{dr}{dt}\right)^2-\frac{r\theta_0^2}{(r-n)(r+n)^3}=0.  \label{16}\\
\end{aligned}
\end{equation}
Together with (\ref{h1}), we have
\begin{align*}
\frac{d}{dt}\left[\frac{r+n}{r-n}\left(\frac{dr}{dt}\right)^2+\frac{\theta_0^2}{r^2-n^2}\right]=0.
\end{align*}
Therefore,
\begin{align*}
\frac{r+n}{r-n}\left(\frac{dr}{dt}\right)^2+\frac{\theta_0^2}{r^2-n^2}=r_1^2.
\end{align*}
Then we obtain
\begin{align*}
\frac{dr}{dt}=\frac{\ep r_1}{r+n}\sqrt{r^2-R_3^2}.
\end{align*}
Thus,
\beQ
\begin{aligned}
\frac{dt}{dr}&=\frac{\ep(r+n)}{r_1\sqrt{r^2-R_3^2}},\\
\frac{d\theta}{dr}&=\frac{d\theta}{dt}\frac{dt}{dr}=\frac{\ep\theta_0}{r_1(r-n)\sqrt{r^2-R_3^2}}.
\end{aligned}
\eeQ
The theorem follows by integrating them from $R_3$ to $r$.\qed
\subsection{Geodesics for constant $r$}\ls
The geodesic equations reduce to
\begin{align}
\frac{d^2\tau}{dt^2}&+\frac{4n^2\cos\theta}{(r+n)^2\sin\theta}\frac{d\tau}{dt}\frac{d\theta}{dt}\nonumber\\
&+\frac{2[4n^3\cos^2\theta-n\sin^2\theta(r+n)^2-2n(r+n)^2\cos^2\theta]}{(r+n)^2\sin\theta}\frac{d\varphi}{dt}\frac{d\theta}{dt}=0,\nonumber\\
\frac{d^2\varphi}{dt^2}&+2\left[-\frac{2n^2\cos\theta}{(r+n)^2\sin\theta}+\frac{\cos\theta}{\sin\theta}\right]\frac{d\varphi}{dt}\frac{d\theta}{dt}
-\frac{2n}{(r+n)^2\sin\theta}\frac{d\tau}{dt}\frac{d\theta}{dt}=0,\nonumber\\
\frac{d^2\theta}{dt^2}&+\left[\frac{4n^2\cos\theta\sin\theta}{(r+n)^2}-\sin\theta\cos\theta\right]\left(\frac{d\varphi}{dt}\right)^2
+\frac{2n\sin\theta}{(r+n)^2}\frac{d\tau}{dt}\frac{d\varphi}{dt}=0,\nonumber\\
&\frac{n(r-n)}{(r+n)^3}\left(\frac{d\tau}{dt}\right)^2+\left[\frac{4n^3\cos^2\theta}{(r+n)^2}+r\sin^2\theta\right]\frac{r-n}{r+n}\left(\frac{d\varphi}{dt}\right)^2\nonumber\\
&+
\frac{r(r-n)}{r+n}\left(\frac{d\theta}{dt}\right)^2
+\frac{4n^2(r-n)\cos\theta}{(r+n)^3}\frac{d\tau}{dt}\frac{d\varphi}{dt}=0.\label{19}
\end{align}
(\ref{19}) implies that
\begin{align*}
\frac{n(r-n)}{(r+n)^3}\left(\frac{d\tau}{dt}\right.&\left.+2n\cos\theta\frac{d\varphi}{dt}\right)^2
+\frac{r(r-n)\sin^2\theta}{r+n}\left(\frac{d\varphi}{dt}\right)^2\\
&+\frac{r(r-n)}{r+n}\left(\frac{d\theta}{dt}\right)^2=0.
\end{align*}
Therefore,
\begin{align*}
\frac{d\tau}{dt}=\frac{d\varphi}{dt}=\frac{d\theta}{dt}=0.
\end{align*}
Thus, $\tau$, $\varphi$, $\theta$ are constants.
\subsection{Geodesics for constant $\varphi$}\ls
The geodesic equations reduce to
\begin{align}
\frac{d^2\tau}{dt^2}&+\frac{4n^2\cos\theta}{(r+n)^2\sin\theta}\frac{d\tau}{dt}\frac{d\theta}{dt}+\frac{2n}{r^2-n^2}\frac{d\tau}{dt}\frac{dr}{dt}=0,\nonumber\\
\frac{d^2\theta}{dt^2}&+\frac{2r}{r^2-n^2}\frac{dr}{dt}\frac{d\theta}{dt}=0,\nonumber\\
\frac{d^2r}{dt^2}&-\frac{n}{r^2-n^2}\left(\frac{dr}{dt}\right)^2-\frac{n(r-n)}{(r+n)^3}\left(\frac{d\tau}{dt}\right)^2-\frac{r(r-n)}{r+n}\left(\frac{d\theta}{dt}\right)^2=0,\nonumber\\
&-\frac{2n}{(r+n)^2\sin\theta}\frac{d\tau}{dt}\frac{d\theta}{dt}=0.\label{26}
\end{align}
(\ref{26}) implies that
\begin{align*}
\frac{d\tau}{dt}=0\,\,\text{or}\,\, \frac{d\theta}{dt}=0.
\end{align*}
It reduces to the cases in Sect. 3.2 or Sect. 3.3 respectively.
\subsection{Geodesics for constant $\theta$}\ls
In this subsection, we solve the geodesic equations for constant $\theta\in(0,\pi)$ on self-dual Taub-NUT metric given by (\ref{ge}). In this case, the geodesic equations reduce to
\begin{align}
\frac{d^2\tau}{dt^2}&+\frac{2n}{r^2-n^2}\frac{d\tau}{dt}\frac{dr}{dt}-\frac{4n\cos\theta}{r+n}\frac{d\varphi}{dt}\frac{dr}{dt}=0,\nonumber\\
\frac{d^2\varphi}{dt^2}&+\frac{2r}{r^2-n^2}\frac{d\varphi}{dt}\frac{dr}{dt}=0,\label{4-2}\\
\frac{d^2r}{dt^2}&-\frac{n}{r^2-n^2}\left(\frac{dr}{dt}\right)^2-\frac{n(r-n)}{(r+n)^3}\left(\frac{d\tau}{dt}\right)^2\nonumber\\
&-\left[\frac{4n^3\cos^2\theta}{(r+n)^2}+r\sin^2\theta\right]\frac{r-n}{r+n}\left(\frac{d\varphi}{dt}\right)^2-\frac{4n^2(r-n)\cos\theta}{(r+n)^3}\frac{d\tau}{dt}\frac{d\varphi}{dt}=0,\label{4-1}\\
&-\left[\sin\theta\cos\theta-\frac{4n^2\cos\theta\sin\theta}{(r+n)^2}\right]\left(\frac{d\varphi}{dt}\right)^2+\frac{2n\sin\theta}{(r+n)^2}\frac{d\tau}{dt}\frac{d\varphi}{dt}=0.\label{20}
\end{align}
(\ref{20}) implies that
\begin{align*}\frac{d\varphi}{dt}=0,
\end{align*}
or \begin{align*}
\left[\frac{4n^2\cos\theta\sin\theta}{(r+n)^2}-\sin\theta\cos\theta\right]\frac{d\varphi}{dt}+\frac{2n\sin\theta}{(r+n)^2}\frac{d\tau}{dt}=0.
\end{align*}
It reduces to Theorem \ref{t2} if $\frac{d\varphi}{dt}=0$.

Now we focus on the case
\begin{align*}
\left[\frac{4n^2\cos\theta}{(r+n)^2}-\cos\theta\right]\frac{d\varphi}{dt}+\frac{2n}{(r+n)^2}\frac{d\tau}{dt}=0.
\end{align*}
Let $r_1, \varphi_0$ be constant and $r_1\neq 0$. Denote
\begin{align}
F(r)=2nr_1^2r^2-\varphi_0^2\cos^2\theta r-(2n^3r_1^2+n\varphi_0^2\cos^2\theta+2n\varphi_0^2\sin^2\theta).\label{F}
\end{align}
The two roots of $F(r)=0$ are given by
\begin{align*}
R_{\pm}=\frac{\varphi_0^2\cos^2\theta}{4nr_1^2}\pm\sqrt{n^2+\frac{\varphi_0^4\cos^4\theta+8n^2\varphi_0^2
r_1^2\cos^2\theta+16n^2\varphi_0^2r_1^2\sin^2\theta}{16n^2r_1^4}},
\end{align*}
respectively.
\begin{lem}\label{l2}
With the above notations, $R_+\geq n$, and $R_+=n$ if and only if $\varphi_0=0$.
\end{lem}
\pf A straightforward computation.\qed
\begin{thm}\label{t5}
The geodesics for metrics of self-dual Taub-NUT type with constant $\theta$ with conditions
\begin{align*}
&\lim\limits_{r \rw R_+} t= t_{1},\,\,\lim\limits_{r \rw R_+} \tau = \tau_{1},\,\,\lim\limits_{r \rw R_+} {(r^2-n^2)\frac{d\varphi}{dt}} = \varphi_{0},\\
&\lim\limits_{r\rw R_+}{\varphi}=\varphi_1,\,\,\lim\limits_{r\rw R_+}{\left[\frac{r+n}{r-n}\left(\frac{dr}{dt}\right)^2+\frac{\varphi_0^2\cos^2\theta}{2n(r-n)}
+\frac{\varphi_0^2\sin^2\theta}{r^2-n^2}\right]}=r_1^2.
\end{align*}
satisfy
\begin{align*}
t(r)=t_1&+\frac{\ep}{\sqrt{2n}}\Bigg[\sqrt{(r-R_+)(r-R_-)}\\
&+\left(\frac{\varphi_0^2\cos^2\theta}{2nr_1^2}+2n\right)\sinh^{-1}\left(\sqrt{\frac{r-R_+}{R_+-R_-}}\right)\Bigg],\\
\varphi(r)=\varphi_1&+\frac{2\ep\sqrt{2n}\varphi_0}{\sqrt{(R_+-n)(n-R_-)}}\arctan\left(\sqrt{\frac{(r-R_+)(n-R_-)}{(r-R_-)(R_+-n)}}\right),\\
\tau(r)=\tau_1&+\frac{\ep\varphi_0\cos\theta}{\sqrt{2n}}\Bigg[\sqrt{(r-R_+)(r-R_-)}\\
&+\left(\frac{\varphi_0^2\cos^2\theta}{2nr_1^2}+6n\right)\sinh^{-1}\left(\sqrt{\frac{r-R_+}{R_+-R_-}}\right)\Bigg].
\end{align*}
\end{thm}
\pf By using the argument in Sect. 3.3, we obtain that (\ref{4-2}) implies
\begin{align*}
\frac{d\varphi}{dt}=\frac{\varphi_0}{r^2-n^2}.
\end{align*}
This gives
\begin{align*}
\frac{d\tau}{dt}=\frac{\varphi_0(r+3n)\cos\theta}{2n(r+n)}.
\end{align*}
Substituting them into (\ref{4-1}), we obtain
\begin{align*}
\frac{d^2r}{dt^2}&-\frac{n}{r^2-n^2}\left(\frac{dr}{dt}\right)^2-\frac{\varphi_0^2r\sin^2\theta}{(r+n)^3 (r-n)}-\frac{\varphi_0^2\cos^2\theta}{4n(r^2-n^2)}=0.
\end{align*}
Then we have
\begin{align*}
\frac{d}{dt}\left[\frac{r+n}{r-n}\left(\frac{dr}{dt}\right)^2+\frac{\varphi_0^2\cos^2\theta}{2n(r-n)}
+\frac{\varphi_0^2\sin^2\theta}{r^2-n^2}\right]=0.
\end{align*}
Thus,
\begin{align*}
\frac{r+n}{r-n}\left(\frac{dr}{dt}\right)^2+\frac{\varphi_0^2\cos^2\theta}{2n(r-n)}
+\frac{\varphi_0^2\sin^2\theta}{r^2-n^2}=r_1^2.
\end{align*}
This implies that
\begin{align*}
\frac{dr}{dt}=\frac{\ep\sqrt{F(r)}}{\sqrt{2n}(r+n)}=\frac{\ep\sqrt{(r-R_+)(r-R_-)}}{\sqrt{2n}(r+n)}.
\end{align*}
Therefore,
\begin{align*}
\frac{dt}{dr}&=\frac{\ep\sqrt{2n}(r+n)}{\sqrt{(r-R_+)(r-R_-)}},\\
\frac{d\varphi}{dr}&=\frac{d\varphi}{dt}\frac{dt}{dr}=\frac{\ep\sqrt{2n}\varphi_0}{r-n}\frac{1}{\sqrt{(r-R_+)(r-R_-)}},\\
\frac{d\tau}{dr}&=\frac{d\tau}{dt}\frac{dt}{dr}=\frac{\ep\varphi_0(r+3n)\cos\theta}{\sqrt{2n}}\frac{1}{\sqrt{(r-R_+)(r-R_-)}}.
\end{align*}
The theorem follows by integrating the above equations from $R_+$ to $r$.\qed
\footnotesize {

\noindent {\bf Acknowledgement} The work is supported by the National Natural Science Foundation of China 12301072.

}

\end{document}